\documentclass[12pt]{amsart}
\pdfoutput=1
\usepackage[T1]{fontenc}
\usepackage[english]{babel}
\usepackage{graphicx}
\usepackage{hyperref}
\usepackage{dsfont}
\usepackage{soul}
\usepackage{enumitem}
\usepackage{tikz-cd}
\graphicspath{{graphics/}}
\usepackage{geometry}
\setlist[enumerate]{labelsep=*, leftmargin=1.5pc,
topsep=1ex plus0.5ex minus0.2ex,
itemsep=1ex plus0.5ex minus0.2ex,
font=\rmfamily,
font=\upshape}
\setlist[itemize]{labelsep=*, leftmargin=1.5pc,
topsep=1ex plus0.5ex minus0.2ex,
itemsep=1ex plus0.5ex minus0.2ex,
font=\rmfamily,
font=\upshape}
\newtheorem{thm}{Theorem}[section]

\newtheorem{lem}[thm]{Lemma}

\newtheorem*{thm*}{Theorem}
\theoremstyle{definition}

\newtheorem{exa}[thm]{Example}

\numberwithin{equation}{section}
\renewcommand{\Re}{\operatorname{Re}}
\renewcommand{\Im}{\operatorname{Im}}
\newcommand{\id}{\mathds{1}}
\newcommand{\ii}{\operatorname{i}}
\newcommand{\tr}{\operatorname{tr}}
\newcommand{\conv}{\operatorname{conv}}
\newcommand{\bC}{\mathbb C}
\newcommand{\bE}{\mathbb E}
\newcommand{\bN}{\mathbb N}
\newcommand{\bP}{\mathbb P}
\newcommand{\bR}{\mathbb R}
\newcommand{\bS}{\mathbb S}

\newcommand{\cC}{\mathcal C}
\newcommand{\cX}{\mathcal X}
\newcommand{\cB}{\mathcal B}
\begin{document}
\selectlanguage{english}
\title{On a theorem by Kippenhahn}
\author{Stephan Weis}
\begin{abstract}   
Kippenhahn discovered a real algebraic plane curve whose convex hull 
is the numerical range of a matrix. The correctness of this theorem 
was called into question when Chien and Nakazato found an example 
where the spatial analogue fails. They showed that the mentioned 
plane curve indeed lies inside the numerical range. We prove the 
easier converse direction of the theorem. Finding higher-dimensional 
generalizations of Kippenhahn's theorem is a challenge in real 
algebraic geometry.
\end{abstract}
%
%
%
%
%
%
\maketitle
\thispagestyle{empty}
\pagestyle{myheadings}
\markleft{\hfill 
On a theorem by Kippenhahn
\hfill}
\markright{\hfill S.~Weis\hfill}
%
%
\section{Introduction}
\label{sec:intro}
\par
Geometric methods play a major role in the theory of the 
{\em numerical range}\footnote{%
We use Dirac's notation \cite{NielsenChuang2010} for the inner product 
of $\bC^d$.}
\[
\{\langle\psi|A\psi\rangle :
|\psi\rangle\in\bC^d, \langle \psi|\psi\rangle=1\}
\]
of a complex $d$-by-$d$ matrix $A$. Already Toeplitz and Hausdorff 
\cite{Toeplitz1918,Hausdorff1919} showed that the numerical range is a 
convex subset of the complex plane $\bC$. Kippenhahn's theorem 
\cite{Kippenhahn1951}\,\footnote{%
The article \cite{Kippenhahn1951} is translated into English
\cite{ZachlinHochstenbach2008}. The name of the beautiful town and 
UNESCO World Heritage Site which is misspelled in the title of 
\cite{ZachlinHochstenbach2008} is {\em Bamberg}.}
affirms that the numerical range is the convex hull of an algebraic curve. 
His theorem has been very influential for decades, see for example 
\cite{Fiedler1981,Keeler-etal1997,JoswigStraub1998,Mirman1998, 
ChienNakazato1999,Bebiano-etal2005}, and it is an early example of the 
young field of convex algebraic geometry 
\cite{Blekherman-etal2013,Netzer2012,SinnSturmfels2015}. His theorem is
fundamental also for geometric approaches to quantum mechanics because 
the numerical range is a linear image of the set of quantum states 
\cite{BerberianOrland1967,Dunkl-etal2011}. Therefore, the numerical range 
is a reduced statistical model of quantum mechanics \cite{Holevo2011}. 
Higher-dimensional linear images of the set of quantum states can represent 
quantum marginals, which are essential for the theory of quantum many-body 
systems \cite{Coleman1963,Erdahl1972}.
\par
Higher-dimensional generalizations of the numerical range were studied in 
operator theory. We introduce them from an algebraic view and consider the 
*-algebra of complex $d$-by-$d$ matrices $M_d$, $d\in\bN$, endowed with the 
Hilbert-Schmidt inner product $\langle A,B\rangle=\tr(A^*B)$ of 
$A,B\in M_d$. We define  
\begin{equation}\label{eq:dm}
\cB=\{\rho\in M_d : \tr(\rho)=1, \rho\succeq 0\},
\end{equation}
where $A\succeq 0$ means that $A\in M_d$ is positive semi-definite. Elements
of $\cB$ are called {\em density matrices} or {\em mixed states}
\cite{BengtssonZyczkowski2006}. The set $\cB$ is a convex and compact set 
called {\em state space} \cite{AlfsenShultz2001} of $M_d$. For $k\in\bN$ 
hermitian matrices $F_1,\ldots,F_k\in M_d$, the convex set 
\[
W=W_{F_1,\ldots,F_k}
:=\{(\langle F_1,A\rangle,\ldots,\langle F_k,\rho\rangle) 
 : \rho\in\cB\}\subset\bR^k,
\]
is known as the {\em joint algebraic numerical range} \cite{Mueller2010} 
of $F_1,\ldots,F_k$, which we call {\em convex support} \cite{Weis2011} in 
analogy with statistics \cite{Barndorff-Nielsen1978}. 
\par
To get in touch with the work \cite{ChienNakazato2010} by Chien and 
Nakazato we define the {\em joint numerical range} of $F_1,\ldots,F_k$ by
\[
\{(\langle\psi|F_1\psi\rangle,\ldots,\langle \psi|F_k\psi\rangle) :
|\psi\rangle\in\bC^d, \langle \psi|\psi\rangle=1\}\subset\bR^k.
\]
For $k=2$, if we identify $\bR^2\cong\bC$, then the joint numerical range is 
the numerical range of $F_1+\ii F_2$, which is convex by the Toeplitz-Hausdorff 
result. For $k=3$ and $d\geq 3$ the joint numerical range is also convex 
\cite{Au-YeungPoon1979}, but there is no simple rule of convexity for 
$k\geq 4$, see for example \cite{LiPoon2000}. However, for all $d,k\in\bN$ 
and hermitian matrices $F_1,\ldots,F_k\in M_d$, the convex hull of the joint 
numerical range is the convex support $W_{F_1,\ldots,F_k}$, see for example 
\cite{Mueller2010,Szymanski-etal2016}. For $k=2$ this means 
\cite{BerberianOrland1967} that the numerical range of $F_1+\ii F_2$ 
is $W_{F_1,F_2}$. 
\par
Algebraic geometry is a natural part of the geometry of $W$ because the 
determinant of a linear matrix polynomial vanishes on the boundary of the 
convex dual of $W$ (see the appendix). More specifically, by reconstructing 
$W$ from this determinant polynomial, which is our declared goal, we enter 
the field of {\em real algebraic geometry} \cite{Bochnak-etal1998}. For 
$u=(u_0,\ldots,u_k)\in\bC^{k+1}$ let 
\begin{equation}\label{eq:linm}
L(u)=L_{F_1,\ldots,F_k}(u):=u_0\id+u_1F_1+\cdots+u_kF_k,
\end{equation}
where $\id\in M_d$ is the identity matrix. The polynomial $\det(L(u))$ 
is homogeneous in $u$ and defines the complex projective variety 
\begin{equation}\label{eq:det-variety}
V=V_{F_1,\ldots,F_k}
:=\{(u_0:u_1:\ldots:u_k)\in\bP^k(\bC)\mid \det(L(u))=0 \}.
\end{equation}
The {\em boundary generating hypersurface} \cite{ChienNakazato2010} of 
$F_1,\ldots,F_k$ is defined by 
\[
V^\vee_\bR
=V^\vee_{\bR,F_1,\ldots,F_k}
=\{(x_1,\ldots,x_k)\in\bR^k \mid (1:x_1:\ldots:x_k)\in V^\vee \},
\]
where $V^\vee=V^\vee_{F_1,\ldots,F_k}\subset\bP^k(\bC)^*$ is the 
{\em projective dual} \cite{Gelfand-etal1994} of $V$. If $k=2$ then 
$V_{F_1,F_2}$ is an {\em algebraic curve}, its projective dual 
$V^\vee_{F_1,F_2}$ is the {\em dual curve} \cite{Fischer2001} of 
$V_{F_1,F_2}$, and the real algebraic set $V^\vee_{\bR,F_1,F_2}$ is 
the {\em boundary generating curve} \cite{Kippenhahn1951} of 
$F_1+\ii F_2$.  
\par
The inclusion $V^\vee_{\bR,F_1,F_2}\subset W_{F_1,F_2}$ is proved in
Corollary~2.5 of \cite{ChienNakazato2010} using Puiseux series. Together 
with the inclusion $W_{F_1,F_2}\subset\conv(V^\vee_{\bR,F_1,F_2})$, 
which we show in Lemma~\ref{lem:ext-on-curve}, Kippenhahn's theorem
\cite{Kippenhahn1951} follows:
\begin{thm}\label{thm:kippenhahn}
Let $d\in\bN$ and $F_1,F_2\in M_d$ be hermitian matrices.
Then $W_{F_1,F_2}=\conv(V^\vee_{\bR,F_1,F_2})$.
\end{thm}
\par
The boundary generating hypersurface $V^\vee_{\bR,F_1,\ldots,F_k}$ becomes 
more interesting for $k\geq 3$. Indeed, an example from \cite{ChienNakazato2010} 
with $k=3$ hermitian $3$-by-$3$ matrices $F_1,F_2,F_3\in M_3$ shows that 
$V^\vee_{\bR,F_1,F_2,F_3}$ can contain lines. This makes the inclusion 
$V^\vee_{\bR,F_1,F_2,F_3} \subset W_{F_1,F_2,F_3}$ impossible because 
$W_{F_1,F_2,F_3}$ is compact. Further, $V^\vee_{\bR,F_1,F_2,F_3}$ is an 
irreducible surface of dimension two and the line has local dimension one. So 
the failure is a phenomenon of real algebraic geometry because an irreducible 
complex algebraic subset of $\bC^k$ has constant local dimension, see 
Section~3.1 of \cite{Bochnak-etal1998}. 
\begin{exa}\label{exa:roman}
The boundary generating surface of the matrices 
\[
F_1:=\tfrac{1}{2}
\left(\begin{smallmatrix}
 0 & 1 & 0 \\
 1 & 0 & 0 \\
 0 & 0 & 0 \\
\end{smallmatrix}\right),
\quad
F_2:=\tfrac{1}{2}
\left(\begin{smallmatrix}
 0 & 0 & 1\\
 0 & 0 & 0\\
 1 & 0 & 0 \\
\end{smallmatrix}\right),
\quad
F_3:=\tfrac{1}{2}
\left(\begin{smallmatrix}
 0 & 0 & 0 \\
 0 & 0 & 1 \\
 0 & 1 & 0 \\
\end{smallmatrix}\right)
\] 
is the {\em Roman surface} 
$\{(x_1,x_2,x_3)\in\bR^3:x_1 x_2 x_3 - x_1^2 x_2^2 - x_1^2 x_3^2 - x_2^2 x_3^2=0\}$,
which contains the three coordinate axes. 
\end{exa}
\par
Example~\ref{exa:roman} is taken from the collection of examples in 
Section~6 of \cite{Szymanski-etal2016}, which suggests that all unbounded 
Zariski closed proper subsets of $V^\vee_{\bR,F_1,F_2,F_3}$ are lines for 
$d=3$. Towards a generalization of Theorem~\ref{thm:kippenhahn} from 
$k=2$ to $d=k=3$, a natural question would be whether the remainder 
of $V^\vee_{\bR,F_1,F_2,F_3}$ without lines is included in 
$W_{F_1,F_2,F_3}$. It could also be helpful to clarify relationships 
between $V^\vee_{\bR,F_1,F_2,F_3}$ and the {\em algebraic boundary} 
\cite{RostalskiSturmfels2013,SinnSturmfels2015} of $W_{F_1,F_2,F_3}$, 
the smallest complex algebraic variety that contains the boundary of 
$W_{F_1,F_2,F_3}$ in the topology of $\bR^3$. 
%
%
\section{Algebraic duality}
\label{sec:alg_duality}
\par
We prove the inclusion of the numerical range into the convex hull 
of the boundary generating curve.
\par
The {\em real part} of a matrix $A\in M_d$ is $\Re(A)=\tfrac{1}{2}(A+A^*)$,
the {\em imaginary part} is $\Im(A)=\tfrac{1}{2\ii}(A-A^*)$. Let 
$F_1,F_2\in M_d$ be hermitian matrices and put $A=F_1+\ii F_2$. Let
\[
A(\theta):=\Re(e^{-\ii\theta}A)
=\cos(\theta)F_1+\sin(\theta)F_2,
\qquad
\theta\in\bR,
\]
and let $\bS(\bC^d):=\{|\psi\rangle\in\bC^d : \langle\psi|\psi\rangle=1 \}$ 
denote the unit sphere. The {\em numerical range map} of $A$ is defined by
$f_A:\bS(\bC^d)\to\bC$,
$|\psi\rangle \mapsto \langle \psi|A\psi\rangle$.
An easy calculation shows 
\begin{equation}\label{eq:NR-map}
f_A(|\psi\rangle)
=e^{\ii\theta}[f_{A(\theta)}(|\psi\rangle)+\ii f_{A'(\theta)}(|\psi\rangle)],
\qquad
|\psi\rangle\in \bS(\bC^d), \theta\in\bR,
\end{equation}
where $A'(\theta)$ denotes derivative with respect to $\theta$
(notice that $A'(\theta)=\Im(e^{-\ii\theta}A)$). An {\em extreme point} $z$ 
of $W$ is a point of $W$ which cannot be written as a proper convex 
combination of other points of $W$, that is if $x,y\in W$, $s\in(0,1)$, and 
$z=(1-s)x+s y$, then $x=y=z$.
\begin{lem}\label{lem:ext-on-curve}
Let $d\in\bN$ and $F_1,F_2\in M_d$ be hermitian matrices. 
Then $W_{F_1,F_2}\subset\conv(V^\vee_{\bR,F_1,F_2})$.
\end{lem}
{\em Proof:} 
As a compact convex set is the convex hull of its extreme points, see for example 
Section~18 of \cite{Rockafellar1970}, it suffices to show that every extreme point 
$z=z_1+\ii z_2$ of $W=W_{F_1,F_2}$ lies on the boundary generating curve 
$V^\vee_\bR=V^\vee_{\bR,F_1,F_2}$. 
\par
Without loss of generality we assume that $z$ lies on the supporting line of 
$W$  which is parallel to the imaginary axis and which meets $W$ on the left
(notice that $1\in\bC$ is an inner normal vector of $W$ at $z$). Under these 
assumptions, Corollaries~2.4 and~2.5 of \cite{SpitkovskyWeis2016} show that 
there exist 
\[
|\psi(\theta)\rangle\in\bS(\bC^d)
\quad\mbox{and}\quad
\lambda(\theta)\in\bR,
\qquad\theta\in\bR,
\]
both analytic in $\theta$, such that 
\begin{equation}\label{eq:z-on-curve}
f_A(|\psi(0)\rangle)=z
\end{equation} 
and for all $\theta\in\bR$
\begin{equation}\label{eq:ev-curve}
A(\theta)|\psi(\theta)\rangle=\lambda(\theta)|\psi(\theta)\rangle.
\end{equation} 
The proof of (\ref{eq:z-on-curve}) and (\ref{eq:ev-curve}) uses perturbation 
theory \cite{Rellich1954} and properties of the support function 
\cite{Schneider2014}, which denotes the scaled distance of the origin from 
the supporting line of $W$ with inner normal vector $u\in\bC$. Equation 
(\ref{eq:ev-curve}) shows
\[
\det[-\lambda(\theta)\id+\cos(\theta)F_1+\sin(\theta)F_2]=0,
\]
so 
\[
u(\theta):=[-\lambda(\theta):\cos(\theta):\sin(\theta)]
\]
lies on the algebraic curve $V=V_{F_1,F_2}$ defined in (\ref{eq:det-variety}).
If $u(\theta)$ is a smooth point of $V$ then the tangent space of $V$ at 
$u(\theta)$ is represented by
\[
x(\theta):=[1:\lambda(\theta)\cos(\theta)-\lambda'(\theta)\sin(\theta):
\lambda(\theta)\sin(\theta)+\lambda'(\theta)\cos(\theta)]\in
\bP^2(\bC)^*.
\]
Since $u(\theta)$ is nowhere locally constant, we get 
$x(0)\in V^\vee=V^\vee_{F_1,F_2}$. Indeed, the dual curve $V^\vee$ is the 
closure of the set of tangent spaces at smooth points of the algebraic curve 
$V$, and $V$ has at most finitely many singular points, see 
Section~3.2 of \cite{Fischer2001}. We observe that
\[
\lambda(0)=f_{A(0)}(|\psi(0)\rangle)
\]
and cite from Lemma~3.2 of \cite{JoswigStraub1998} the equation
\[
\lambda'(0)=f_{A'(0)}(|\psi(0)\rangle).
\]
So (\ref{eq:NR-map}) and (\ref{eq:z-on-curve}) show 
\[
x(0)
=[1:\lambda(0):\lambda'(0)]
=[1:\Re(z):\Im(z)].
\]
This shows $z\in V^\vee_{\bR,F_1,F_2}$ and finishes the proof.
\hspace*{\fill}$\square$\\
\par
We remark that the representation of {\em extreme points} of $W$, 
provided in \cite{SpitkovskyWeis2016}, can be replaced in the proof of 
Lemma~\ref{lem:ext-on-curve} with the simpler and more well-known 
representation of {\em exposed points} of $W$. In fact, $W$ is closed
and Straszewicz's theorem \cite{Schneider2014} shows that the extreme 
points of $W$ are covered by the closure of the set of exposed points 
of $W$.
%
%
\section{Appendix: Convex duality}
\label{sec:convex_duality}
\par 
We discuss (partial) dualities for convex sets constructed from a self-dual 
convex cone. As a special case we explain the convex duality between convex 
support sets and spectrahedra. Some definitions from earlier section are 
renewed.
\par
Let $\bE$ be a Euclidean vector space with scalar product 
$\langle x,y\rangle$ of $x,y\in\bE$ and norm $\|x\|=\sqrt{\langle x,x\rangle}$. 
We use the scalar product $\langle u,x\rangle=u_1x_1+\cdots+u_kx_k$ of 
$u,x\in\bR^k$, where $u=(u_1,\ldots,u_k)$ and $x=(x_1,\ldots,x_k)$. For any 
subset $\cX\subset\bE$, the {\em (convex) dual} of $\cX$ is defined by
\[
\cX^*:=\{A\in\bE : \forall B\in\cX, 1+\langle A,B\rangle\geq 0\}.
\]
Recall that $\cX^*$ is a closed convex subset of $\bE$ including the origin. 
We assume that $\cC\subset\bE$ is a {\em self-dual convex cone}, that is
\[
\cC=\{A\in\bE : \forall B\in \cC, \langle A,B\rangle\geq 0\}
\]
(then $\cC^*=\cC$ holds and $\cC$ has non-empty interior). For $F_0\in\bE$ 
we define
\[
\cB:=\{A\in\cC : \langle F_0,A \rangle =1 \}.
\]
For $k\in\bN$ and $F_1,\ldots,F_k\in\bE$ we consider the linear 
$\bE$-polynomial 
\[
L(u):=F_0+u_1F_1+\cdots+u_kF_k,
\qquad u\in\bR^k.
\]
Let
\begin{equation}\label{eq:sectra}
S=S_{F_1,\ldots,F_k}:=\{u\in\bR^k : L(u)\in\cC\} \subset\bR^k
\end{equation}
and 
\begin{equation}\label{eq:conv_supp}
W=W_{F_1,\ldots,F_k}
:=\{(\langle F_1,A\rangle,\ldots,\langle F_k,A\rangle) : A\in\cB\}\subset\bR^k.
\end{equation}
Observe that $S$ is closed and convex. Since $\cB\subset\cC$, we have
\begin{align}\label{eq:duality01}
S &\subset \{u\in\bR^k : \forall A\in\cB, \langle L(u),A\rangle\geq 0 \}\\\nonumber
 &= \{u\in\bR^k : \forall A\in\cB, 
 1+u_1\langle F_1,A\rangle+\cdots+u_k\langle F_k,A\rangle \geq 0 \}\\\nonumber
 &= \{u\in\bR^k : \forall x\in W, 
 1+\langle u,x\rangle \geq 0 \}\\\nonumber
 &= W^*.\nonumber
\end{align}
We remark that (\ref{eq:duality01}) and general properties of convex 
duals \cite{Rockafellar1970} show
\begin{equation}\label{eq:weak}
\conv(0,W)\subset S^*.
\end{equation}
The inclusion (\ref{eq:weak}) follows also, analogous to the last inequality 
on page 45 of \cite{RamanaGoldman1995}, from weak duality of linear 
programming. By this we mean that for all $x\in\bR^k$ 
\[
\sup\{\langle x,-u\rangle : u\in S\}
\leq 
\inf\{\langle F_0,A\rangle : x_i=\langle F_i,A\rangle, 
i=1,\ldots,k, A\in\cC\}.
\]
Thereby, in (\ref{eq:duality01}) and (\ref{eq:weak}) the vector $F_0\in\bE$ 
may be arbitrary, $\cB$ may be unbounded, and $W$ and $\conv(0,W)$ may be 
neither closed nor bounded.
\par
To return to the subject of convex support sets, we assume that $F_0$ is 
an interior point of $\cC$. Then $\cB$ is a non-empty compact and convex 
subset of $\cC$. It is easy to show\footnote{%
That $\cB$ is a base of $\cC$ means that for every $A\in\cC$ there is 
$B\in\cB$ and $\lambda\geq 0$ such that $A=\lambda B$.
See Section~2.4 of \cite{Weis2011} for a proof.}
that $\cB$ is a base of $\cC$, so (\ref{eq:duality01}) gives
\begin{equation}\label{eq:duality02}
S = W^*
\qquad \mbox{if $F_0$ is an interior point of $\cC$.}
\end{equation}
We remark that since $W$ is bounded, the origin is an interior point of $S$, 
but $S$ may be unbounded. Since $W$ is closed, one has $W=S^*$ if $W$ 
contains the origin \cite{Rockafellar1970}. 
\par
In the context of the numerical range, the equation (\ref{eq:duality02}) is 
well-known \cite{ChienNakazato2010,HeltonSpitkovsky2012}. Let $\bE\subset M_d$ 
denote the subspace of hermitian matrices. The cone $\cC\subset\bE$ of 
positive semi-definite matrices is self-dual \cite{BermanBen-Israel1973}. If 
$F_0=\id$ then $\cB$ is the set of density matrices (\ref{eq:dm}) and 
$W=W_{F_1,\ldots,F_k}$ is the convex support of $F_1,\ldots,F_k$. The 
polynomial $L$ is the linear matrix polynomial (\ref{eq:linm}) with $u_0=1$
and
\[
S=\{u\in\bR^k : L(u)\succeq 0\} \subset\bR^k
\] 
is known as a {\em spectrahedron} \cite{RamanaGoldman1995}. Clearly, the 
determinant of $L(u)$ is strictly positive in the interior of $S$ and 
vanishes on the boundary of $S$. 
%
%
%
\vspace{2\baselineskip}
\par\normalsize
{\par\noindent\footnotesize
{\em Acknowledgements.} 
The author thanks Ilya Spitkovsky for feedback to a draft
of the manuscript.}
%
%
%
%
\bibliographystyle{plain}

%
%
\vspace{1cm}
\parbox{12cm}{%
Stephan Weis\\
e-mail: {\tt maths@stephan-weis.info}\\[.5\baselineskip]
Centre for Quantum Information and Communication \\
Université libre de Bruxelles \\
50 av.~F.D. Roosevelt \\
1050 Bruxelles, Belgium}
\end{document}